\documentclass[11pt]{article}
\usepackage{amscd,amsmath, amssymb, fancyhdr, url}
\usepackage{dutchcal}
\usepackage{ stackengine, scalerel, mathtools, bm}
\usepackage{fourier}
\usepackage{libertine}
\usepackage[backref=page]{hyperref}
\renewcommand*{\backref}[1]{}
\renewcommand*{\backrefalt}[4]{%
    \ifcase #1 (Not cited.)%
    \or        (Cited on page~#2.)%
    \else      (Cited on pages~#2.)%
    \fi}

\newcommand{\version}{version 2.0,\ \ July 25, 2017}
\newcommand{\f}{\varphi}

\newcommand{\Lie}{\operatorname{Lie}}

\newcommand{\RR}{\mathbb{R}}

\numberwithin{equation}{section}

\def\eqref#1{(\ref{#1})}

\newcommand{\ra}{{\:\longrightarrow\:}}
\newcommand{\Z}{{\mathbb Z}}
\newcommand{\C}{{\mathbb C}}
\newcommand{\R}{{\mathbb R}}

\newcommand{\6}{\partial}
\def\1{\sqrt{-1}\:}
\newcommand{\restrict}[1]{{\left|_{{\phantom{|}\!\!}_{#1}}\right.}}
\newcommand{\cntrct}                
{\hspace{2pt}\raisebox{1pt}{\text{$\lrcorner$}}\hspace{2pt}}
\newcommand{\arrow}{{\:\longrightarrow\:}}

\newcommand{\calo}{{\cal O}}

\renewcommand{\bar}{\overline}
\renewcommand{\phi}{\varphi}
\renewcommand{\epsilon}{\varepsilon}

\renewcommand{\leq}{\leqslant}

\newcommand{\im}{\operatorname{im}}

\newcommand{\const}{\operatorname{\text{\sf const}}}

\newcommand{\Aut}{\operatorname{Aut}}

\newcounter{Mycounter}[section]
\newcounter{lemma}[section]
\newcounter{claim}[section]
\renewcommand{\theclaim}{{Claim \thesection.\arabic{claim}}}
\newcommand{\claim}{%
     \setcounter{claim}{\value{Mycounter}}
     \refstepcounter{claim}
     \stepcounter{Mycounter}
     {\noindent \bf \theclaim:\ }}

\newcounter{sublemma}[section]
\newcounter{corollary}[section]
\newcounter{theorem}[section]
\renewcommand{\thetheorem}{{Theorem \thesection.\arabic{theorem}}}
\newcommand{\theorem}{%
     \setcounter{theorem}{\value{Mycounter}}
     \refstepcounter{theorem}
     \stepcounter{Mycounter}
     {\noindent \bf \thetheorem:\ }}

\newcounter{conjecture}[section]
\newcounter{proposition}[section]
\renewcommand{\theproposition} {{Proposition \thesection.\arabic{proposition}}}
\newcommand{\proposition}{%
     \setcounter{proposition}{\value{Mycounter}}
     \refstepcounter{proposition}
     \stepcounter{Mycounter}
     {\noindent \bf \theproposition:\ }}

\newcounter{definition}[section]
\renewcommand{\thedefinition} {{Definition~\thesection.\arabic{definition}}}
\newcommand{\definition}{%
     \setcounter{definition}{\value{Mycounter}}
     \refstepcounter{definition}
     \stepcounter{Mycounter}
     {\noindent \bf \thedefinition:\ }}

\newcounter{example}[section]
\renewcommand{\theexample}{{Example \thesection.\arabic{example}}}
\newcommand{\example}{%
     \setcounter{example}{\value{Mycounter}}
     \refstepcounter{example}
     \stepcounter{Mycounter}
     {\noindent \bf \theexample:\ }}

\newcounter{remark}[section]
\renewcommand{\theremark}{{Remark \thesection.\arabic{remark}}}
\newcommand{\remark}{%
     \setcounter{remark}{\value{Mycounter}}
     \refstepcounter{remark}
     \stepcounter{Mycounter}
     {\noindent \bf \theremark:\ }}

\newcounter{problem}[section]
\newcounter{question}[section]
\makeatletter

\@addtoreset{equation}{section}
\@addtoreset{footnote}{section}
\makeatother

\def\blacksquare{\hbox{\vrule width 5pt height 5pt depth 0pt}}
\def\endproof{\blacksquare}

\newcommand{\proof}{{\bf Proof: \ }}

\begin{document}

\begin{center}
{\LARGE\bf  Positivity of LCK potential}\\[5mm]
{\large
Liviu Ornea,\footnote{Partially supported by a grant of Ministry of Research and Innovation, CNCS - UEFISCDI, project number PN-III-P4-ID-PCE-2016-0065, within PNCDI III.}
 and 
Misha
Verbitsky\footnote{Misha Verbitsky is partially supported by the 
Russian Academic Excellence Project '5-100'.\\[1mm]
\noindent{\bf Keywords:} locally conformally K\"ahler,
potential, plurisubharmonic,  
holomorphicaly convex, regularized maximum.

\noindent {\bf 2010 Mathematics Subject Classification:} {53C55, 32E05, 32E10.}
}\\[4mm]
}

\end{center}

{\small
\hspace{0.15\linewidth}
\begin{minipage}[t]{0.7\linewidth}
{\bf Abstract} \\ 
Let $M$ be a complex manifold 
and $L$ an oriented real line bundle
on $M$ equipped with a flat connection.
An LCK (``locally conformally K\"ahler'')
form is a closed, positive (1,1)-form taking
values in $L$, and an LCK manifold is one which admits
an LCK form. Locally, any LCK form is expressed
as an $L$-valued pluri-Laplacian of a function
called LCK potential. We consider a manifold
$M$ with an LCK form admitting an LCK potential
(globally on $M$), and prove that $M$ admits a 
positive LCK potential. Then $M$ admits 
a holomorphic embedding to a Hopf manifold,
as shown in \cite{ov_ma_10}.
\end{minipage}
}

\tableofcontents


\section{Introduction}
\label{_Intro_Section_}

There are several equivalent ways to define 
locally conformally K\"ahler (LCK) manifolds,
but the most appropriate for the present paper 
is the following. Let $M$ be a complex manifold 
and $L$ an oriented real line bundle
on $M$ equipped with a flat connection.
An LCK form is a closed, positive (1,1)-form taking
values in $L$, and an LCK manifold is one which admits
an LCK form. For a more explicit and detailed 
exposition, see Subsection \ref{_LCK-intro_Subsection_1_}.

Locally, a K\"ahler form is always equal to $dd^c\phi$
for some function $\phi$ which is called {\bf K\"ahler
potential}. This statement has an analogue for LCK
manifolds. Denote by $d_\theta$ the de Rham differential
on $L$-valued differential forms (\ref{dtheta}) and let 
$d_\theta^c:= I d_\theta I^{-1}$ be its complex conjugate.
Then locally any LCK form is expressed as $d_\theta d_\theta^c\phi$
where $\phi$ is called {\bf LCK potential}.

On a compact complex manifold,
any plurisubharmonic function is constant, and
this means that such a manifold cannot have a global
K\"ahler potential. However, an LCK manifold
might have a global LCK potential -- this 
was first observed in \cite{_Verbitsky_vanishing_} and \cite{ov_ma_10}
and much used since then (\cite{goto}, \cite{ov_jgp_09},
\cite{_OV:_Topology_imrn_}, \cite{ov_imrn_12},  \cite{O} and so on).
In the first papers on this subject, the LCK potential was assumed
to be positive, but then we realized that the
existence of a potential is a cohomological
notion, see \ref{van_bc}, and the focus was shifted on the study
of the corresponding cohomology group.

LCK manifolds with potential enjoy several
properties which make this notion quite useful.
First, unlike the LCK manifolds (and like the K\"ahler
manifolds) the class of LCK manifolds with potential
is ``deformationally stable'': a small deformation
of an LCK manifold with potential is again LCK with 
potential. Also, any LCK manifold $M$ with potential,
$\dim_\C M >2$,
can be holomorphically embedded in a Hopf
manifold $\C^n \backslash \{0\}/\Z$, and, conversely, 
any complex submanifold of a Hopf manifold is LCK
with potential.

To reconcile the cohomological and geometrical
approach, we need to prove that any 
LCK manifold with the LCK-form $\omega$ 
in the image of $d_\theta d^c_\theta$ admits another
LCK potential which is positive, and this
is the aim of the present paper.

The obvious solution which one would use
in K\"ahler case (adding a constant) does not
work, because the operator $d_\theta d^c_\theta$ 
does not vanish on constants. However, we were able
to find a positive function $h$ such that  
$d_\theta d^c_\theta (h)$ is non-negative; adding
$C d_\theta d^c_\theta (h)$ to $\omega$, $C\gg 0$, gives 
us an LCK form with positive potential in
the same cohomology class.

For a while we tried to prove that any
LCK potential for compact LCK manifolds is positive,
but this claim was wrong; we are grateful to Victor Vuletescu
who disabused us of this fallacy.

Vuletescu's example is the following.
Take a Hopf manifold $\C^n \backslash \{0\}/\Z$
where $\Z$ acts by multiplication by $\lambda >1$
and let $L$ be the local system with the same
monodromy. Then the usual flat K\"ahler form
on $\C^n$ can be considered as a closed
Hermitian form with values in the bundle $L^2$.
Its LCK potential is a function $m(z):= |z|^2$.
Any quadratic polynomial on $\C^n$ gives a holomorphic
section of $L^2\otimes_\R \C$; let $v$ be its real part.
Then $d_\theta d^c_\theta (v)=0$, because $v$ is
the real part of a holomorphic section of $L^2$,
and $d_\theta d^c_\theta (m + Ah)= d_\theta d^c_\theta(m)$
is the LCK form on $M$, for any real constant
$A$. However,
for $A$ large enough, the LCK potential $m+Ah$
is non-positive. 

\hfill

Subsections \ref{_LCK-intro_Subsection_1_}  and
\ref{_LCK-intro_Subsection_2_} are devoted to presenting
the precise definitions and results that will be further
used, and to stating the main result of the paper. Section
\ref{_complex_basic_} recalls a classical Remmert Reduction
Theorem. Section \ref{_gluing_}  extends
Demailly's technique of gluing K\"ahler metrics to the LCK
setting. In Section \ref{_nega_} we begin the proof of our
main result, showing that an LCK potential cannot be
strictly negative,  while in Section \ref{_last_} we
prove that an LCK potential $\phi$ which is positive somewhere can
be glued with another one which is positive on the
set where $\phi <0$ to obtain an LCK
potential which is strictly positive.


\section{Introduction to LCK geometry}


\subsection{LCK manifolds}
\label{_LCK-intro_Subsection_1_}

\definition A complex manifold $(M,I)$ of complex dimension at least 2, is called {\bf locally conformally K\"ahler} (LCK)
if it admits a  Hermitian metric $g$ and a closed 1-form $\theta$,  called {\bf the Lee
form}, such that the  
fundamental 2-form $\omega(\cdot,\cdot):=g(\cdot, I\cdot)$
satisfies the integrability condition
\begin{equation}\label{deflck}
d\omega=\theta\wedge\omega,\quad d\theta=0.
\end{equation}

\hfill

\remark As shown by Vaisman (see
\cite[Theorem 2.1]{_Dragomir_Ornea_}),
a compact locally conformally K\"ahler manifold
with $\theta$ non-exact cannot admit
any K\"ahler form. On the other hand, an LCK form
can be made into a K\"ahler one whenever $\theta$
is exact; indeed, if $\theta= df$,
the form $e^{-f}\omega$ is K\"ahler.
Such a manifold is sometimes called
``globally conformally K\"ahler''.
Throughout this paper,
we shall tacitly exclude this case and assume that $\theta$ is non-exact
for any closed LCK manifold we consider.

\hfill

An equivalent definition is given as follows (see \emph{e.g.} \cite{_Dragomir_Ornea_}):

\hfill

\proposition A complex manifold $(M,I)$, of complex dimension at least 2, is LCK if and only if it admits 
 a  covering $\pi:\tilde M\ra M$ endowed with a K\"ahler
 form $\tilde\omega$ 
with the deck group $\Aut_M(\tilde M)$ acting 
on $(\tilde M, \tilde \omega)$ by holomorphic homotheties. 
 
 \hfil
 
 Hence, if $\tau\in \Aut_M(\tilde M)$, then $\tau^*\tilde\omega=c_\tau\cdot\tilde\omega$, where $c_\tau\in \R^{>0}$ is the scale factor. This defines a character 
\begin{equation}\label{chi}
\chi:\Aut_M(\tilde M)\longrightarrow \R^{>0},\quad \chi(\tau)=c_\tau.
\end{equation}

\definition Differential forms $\tilde\eta$ on $\tilde M$ which satisfy $\gamma^* \eta= c_\gamma \tilde\eta$
for any deck transform map $\gamma\in\Aut_M(\tilde M)$, where $c_\gamma:=\chi(\gamma)\neq 1$,  
are called {\bf automorphic}. In particular, $\tilde\omega$ is automorphic.

\hfill

We shall denote by $L$ the flat line bundle on $M$
associated to this character (it is precisely the {\bf
bundle of densities of weight 1} from conformal
geometry). We fix a  trivialization $v$ of $L$
 such that $\theta$ is a connection form in $L$.
The complexification of $L$, 
considered as a holomorphic line bundle,
will be denoted with the same letter (it
will be clear from the context  which one we refer to). 
This holomorphic line bundle is equipped with a natural Hermitian
metric which is constant on $v$. We
shall call $L$ simply the {\bf weight bundle},
which is a standard term in conformal geometry.

\hfill

\definition The rank of the abelian subgroup $\im(\chi)$ in $\R^{>0}$ is called the {\bf LCK rank} of $M$. It equals the rank of the monodromy group of the bundle $(L, \theta)$.

\hfill

\remark\label{dtheta} The cohomology $H^*(M,L)$ of the local system $L$ is isomorphic with the cohomology of the complex $(\Omega^*(M),d_\theta)$, where   $d_\theta:=d-\theta\wedge$, and is called {\bf Morse Novikov cohomology}. 

\hfill

\example The known examples of LCK manifolds include: Hopf
manifolds (\emph{e.g.} \cite{ov_pams_16}), Oeljeklaus-Toma
manifolds with $t=1$ (\cite{ot}), almost all compact
complex surfaces  (\cite{go, bel, bru}). Such examples
fully justify the interest in LCK geometry.
At the moment, there exists only one example
of non-K\"ahler compact surfaces which
do not admit LCK metric, which is one
of the three Inoue surfaces. If the famous
spherical shell conjecture is true, the rest of non-K\"ahler
surfaces are LCK (\cite{_OV:_embedding_new_proof_}).

\subsection{LCK manifolds with potential}
\label{_LCK-intro_Subsection_2_}

In \cite{ov_ma_10} we introduced the folowing subclass of
LCK manifolds, which, as we proved, share essential features
of K\"ahler manifolds: stability at small deformations and
the Kodaira-type embedding theorem, providing a
holomorphic embedding into a Hopf manifold.

\hfill

\definition \label{_LCK_w_p_orig_Definition_}
 An {\bf LCK manifold with  potential} is a manifold
which admits a K\"ahler covering $(\tilde M, \tilde \omega)$ and a 
smooth function $\tilde\phi:\,\tilde M \rightarrow \RR^{>0}$ 
(the {\bf LCK potential}) satisfying the following conditions:
\begin{description}
\item(i)  $\tilde\phi$ is proper, \emph{i.e.} its level sets are compact;
\item(ii) $\tilde\phi$ is automorphic, \emph{i.e.}  
$\tau^* (\phi)=\chi(\tau) \phi$, for all $\tau\in \Aut_M(\tilde M)$.
\item(iii) $\tilde\phi$ is a  K\"ahler potential, \emph{i.e.} 
$dd^c\tilde\phi = \tilde \omega$.
\end{description}
For the geometric interpretation of these conditions, see
Section \ref{_Intro_Section_}.

\hfill

\example {\bf Vaisman manifolds} are LCK manifolds
$(M,g,\theta)$ with paralel Lee form (with respect to the
Levi Civita connection of the LCK metric). The
$\pi^*g$-squared norm of $\pi^*\theta$ is a positive,
automorphic potential on the universal cover $\pi:\tilde
M\ra M$ (\cite{_Verbitsky_vanishing_}). In particular,
diagonal Hopf manifolds are LCK with potential and, by
stability at small deformations, non-diagonal Hopf
manifolds (which are non-Vaisman) are LCK with potential
too, \cite{ov_ma_10}.

On the other hand, Inoue surfaces (\cite{O}), and the
LCK Oeljeklaus-Toma manifolds (\cite{ov_hopf_surf}) are
not LCK  with potential, and hence this subclass is
strict.

\hfill

\remark\label{van_bc} The existence of a LCK potential
immediately implies the vanishing of the class $[\omega]$
in the Bott-Chern cohomology of $M$ with values in $L$,
see \cite{ov_jgp_09}. 

\hfill 

The meaning of the properness condition (i) in the definition is explained by the following equivalence:

\hfill

\proposition (\cite{ov_jgp_16}) \label{_proper_virt_cy_Claim_}
Let $M$ be a compact manifold, 
$\tilde M$ a covering, and 
$\tilde\phi:\,\tilde M \rightarrow \R^{>0}$ 
an automorphic function.
Then $\tilde\phi$ is proper if and only if the deck transform group $\Aut_M(\tilde M)$ of
$\tilde M$ is virtually cyclic (\emph{i.e.} it contains $\Z$ as a finite index subgroup). 

In particular, a compact LCK manifold with potential has LCK rank 1 if and only if the automorphic potential is proper. 

\hfill

\remark\label{pos_neg} Examining the proof, one can see that it equally  works for $\phi<0$: what is important is that $\phi$ does not pass through zero.

\hfill

\ref{_LCK_w_p_orig_Definition_} can be reformulated
to avoid the K\"ahler cover.

\hfill

\proposition (\cite{ov_jgp_16})\label{_LCK_pot_MN_Definition_}
Let $(M,\omega,\theta)$ be an LCK manifold.
Then $M$ is  { LCK  with  potential}
if and only if there exists a positive
function $\phi\in {C}^\infty (M)$ satisfying 
\begin{equation}\label{pot_on_M}
d_\theta d^c_\theta(\phi)=\omega,
\end{equation}
 where $d^c_\theta =
I d_\theta I^{-1}$.

Explicitly, if $\pi^*\theta=d\f$ on $\tilde M$, then
the K\"ahler potential on $M$ is given by
$e^{-f}\cdot\pi^*\phi$.

\hfill

\remark Note that the LCK potential $\phi$ above is
defined on the manifold $M$, which is often compact. It is not a
K\"ahler potential in the usual sense.

\hfill

Automorphic potentials can be
approximated by proper ones, and hence the properness condition in the definition is not essential as long as one is only interested in complex and differential properties (and not in metric ones). 

\hfill

\proposition (\cite{ov_jgp_16})\label{_LCK_appro_Claim_}
Let $(M,\omega, \theta)$ be an LCK manifold, and
$\phi\in {C}^\infty (M)$ a function satisfying 
$d_\theta d^c_\theta(\phi)=\omega$.
Then $M$ admits an LCK structure $(\omega', \theta')$
of LCK rank 1, approximating $(\omega, \theta)$ in ${C}^\infty$-topology.

\hfill

This paper is instead concerned with the positivity of the potential. 
The main result of this paper is the following theorem.

\hfill

\theorem\label{_main_intro_Theorem_}
Let $M$ be an LCK manifold with a K\"ahler covering 
admitting an automorphic K\"ahler potential.
Then $M$ also admits an LCK metric with a
positive automorphic potential.

\hfill

This theorem, which is proven in Section \ref{_last_},
has the following useful corollary (\cite{ov_ma_10}).

\hfill

\theorem
Let $M$ be an LCK manifold manifold with a K\"ahler covering 
admitting an automorphic K\"ahler potential.
Then $M$ admits a holomorphic embedding to a Hopf
manifold.
\endproof

\hfill

The main tool in the proof will be the gluing of LCK
metrics (see Section \ref{_gluing_}) which is based on
Demailly's regularized maximum of two functions.

\hfill

\remark \ref{_main_intro_Theorem_} fills  a gap in the proofs  of the following previous results of ours:

\begin{itemize}
\item \cite[Theorem 1.4]{ov_jgp_09}, where we claimed  the
  existence of an LCK potential and, in fact, we only
  proved the vanishing of $[\tilde\omega]$ in the Bott-Chern
  group $H^{1,1}_{BC}(M,L)$, \emph{i.e.} the existence
  of an automorphic potential which was not necessarily
  positive;
\item \cite[Theorem 2.3]{ov_imrn_12}, where this result
  was used to embed an LCK manifold admitting a
  holomorphic circle action which is not conformal
to isometry to a Hopf manifold. It is proven by taking
an average of the LCK form with respect to the circle
action and noticing that its Bott-Chern class vanishes.
\end{itemize}
Now the results are true as stated.


\section{Remmert Reduction Theorem}\label{_complex_basic_}


 For the sake of completion, we quote here, without
 proof, a classical result we shall use further on.
 
 \hfill

 \theorem (Remmert reduction, \cite{rem})\label{rem_red}\\
Let $X$ be a holomorphically convex
space. Then there exist a Stein space $Y$ and a proper, surjective, holomorphic map
$f : X \ra Y$ such that
\begin{enumerate}
\item[(i)] $f_*\calo_X = \calo_Y$.

Moreover, the fact that $Y$ is Stein and (i) imply:
\item[(ii)] $f$ has connected fibers.
\item[(iii)]  The map $f^* : \calo_Y (Y)\ra \calo_X (X)$ is an isomorphism.
\item[(iv)] The pair $(f, Y)$ is unique up to biholomorphism, \emph{i.e.} for any other pair $(f_0, Y_0)$
with $Y_0$ Stein and property 1., there exists a biholomorphism $g : Y \ra Y_0$ such that
$f_0 = g \circ f$. 
\end{enumerate}



\section{Gluing K\"ahler forms and LCK forms}\label{_gluing_}


\subsection{Regularized maximum of ${d_\theta
  d^c_\theta}$-plurisubharmonic functions}

In \cite{_Demailly_1982_}, the notion
of a {\em regularized maximum} of two functions was
defined as follows.  

\hfill

\definition\label{_regu_max_Definition_}
(\cite{_Demailly_1982_}) Choose $\epsilon >0$, and let 
$\max_\epsilon:\; \R^2\arrow \R$ be a smooth, convex
function, monotonous in both variables,
which satisfies $\max_\epsilon(x, y) = \max(x,y)$
whenever $|x-y|>\epsilon$. Then $\max_\epsilon$ is called
{\bf a regularized maximum}. 

\hfill

\theorem (\cite{_Demailly_1982_}) \label{reg_max} The regularized maximum
of two plurisubharmonic functions is
again plurisubharmonic.

\hfill

\claim\label{_d_theta_reg_max_Claim_}
Let $\theta$ be a closed form on a complex manifold,
and $\phi, \psi$ two $d_\theta d_\theta^c$-plurisubharmonic
functions. Then $\max_\epsilon(\phi, \psi)$ is also
$d_\theta d_\theta^c$-plurisubharmonic; it is
strictly $d_\theta d_\theta^c$-plurisubharmonic if
$\phi$ and $\psi$ are strictly $d_\theta d_\theta^c$-plurisubharmonic.

\hfill

{\bf Proof:} Since this result is local,  
we may always assume that $\theta= d\rho$ for some positive function $\rho$.
Then $d_\theta(\eta) = e^{-\rho}d(e^\rho\eta)$
and 
\[ 
d_\theta d_\theta^c(\max_\epsilon(\phi, \psi))=e^{-\rho}
dd^c (\max_{e^{\rho\epsilon}}(e^\rho \phi,e^\rho \psi).
\]
Since $e^\rho \phi$, $e^\rho \psi$ are plurisubharmonic, the form
$dd^c (\max_\epsilon(e^\rho \phi,e^\rho  \psi))$ is positive.
\endproof

\subsection{Gluing of LCK potentials}

The following procedure is well known; it was 
much used by J.-P. Demailly (see \emph{e.g.}  
\cite{_Demailly_Paun_}), and, in LCK context,
in our paper \cite{_OV:emb_Sasakian_}.
We call it ``Gluing of K\"ahler metrics''.

\hfill

\proposition \label{_gluing_Kahler_Proposition_}
(gluing of K\"ahler metrics)\\
Let $(M, \omega)$ be a K\"ahler manifold, and 
$D\subset M$ a submanifold of the same dimension with smooth
compact boundary such that $\omega = dd^c\phi$ in
a smooth neighbourhood of $D$, with $\phi$ a plurisubharmonic
function. Let $\psi$ be another plurisusubharmonic function
with $\psi=\phi$ on $\6 D$. Consider a vector field 
$X\in TM\restrict{\6 D}$ which is normal and
outward-pointing everywhere in $\6 D$, and
let $\Lie_X$ denote the derivative of a function
along $X$. Assume that $\Lie_X \psi < \Lie_X \phi$ everywhere
on $\6 D$. Let $D_-$ be an open subset of $D$
which does not intersect a neighbourhood $U$ of $\6 D$,
and $D_+$ an open subset of $M \backslash D$ which
does not intersect $U$. Then there exists a K\"ahler form
$\omega_1$ which is equal to $\omega$ on $D_+$
and to $dd^c \psi$ on $D_-$.

\hfill

\proof
Consider the function $\max_\epsilon(\psi, \phi)$
defined as in \ref{_regu_max_Definition_}, 
where $\psi, \phi$ are defiend as above. Since
$\Lie_X \psi < \Lie_X \phi$, the maximum
$\max(\psi, \phi)$ is equal to $\psi$ in
$D_-$ near $\6 D$ and equal to $\phi$
$D_+$ near $\6 D$. We choose $D_+, D_-$
sufficiently big in such a way that
$\psi > \phi$ in a neighbourhood of the
boundary  $\6 D_-$ and $\psi < \phi$ in a neighbourhood of the
boundary $\6 D_+$. This gives $\max(\phi, \psi)=\phi$
on $D_+$ and $\max(\phi, \psi)=\psi$ on $D_-$.

Choosing $\epsilon$
sufficiently small, the same would hold for
for the regularized maximum $\max_\epsilon(\psi, \phi)$.
Now we can extend $dd^c \max_\epsilon(\psi, \phi)$
to $D_-$ as $dd^c \psi$ and to $D_+$ as $\omega$.
\endproof

\hfill

Replacing $d, d^c$ by $d_\theta$ and $d^c_\theta$ and
using the regularized maximum of 
$d_\theta d^c_\theta$-plurisubharmonic functions as 
in \ref{_d_theta_reg_max_Claim_}, we 
obtain the following LCK-version of this result;
the proof is the same after we replace $d, d^c$ 
by $d_\theta$ and $d^c_\theta$ (note that below $\phi, \psi$ denote functions on $M$, and not on its K\"ahler covering).

\hfill

\proposition \label{_gluing_LCK_Proposition_}
(gluing of LCK metrics)\\
Let $(M, \theta, \omega)$ be an LCK manifold, and 
$D\subset M$ a submanifold of the same dimension with smooth
compact boundary such that $\omega = d_\theta d^c_\theta \phi$ in
a smooth neighbourhood of $D$, with $\phi$ a  
$d_\theta d^c_\theta$-plurisubharmonic
function. Let $\psi$ be another 
$d_\theta d^c_\theta$-plurisusubharmonic function
with $\psi=\phi$ on $\6 D$. Consider a vector field 
$X\in TM\restrict{\6 D}$ which is normal and
outward-pointing everywhere in $\6 D$, and
let $\Lie_X$ denote the derivative of a function
along $X$. Assume that $\Lie_X \psi < \Lie_X \phi$ everywhere
on $\6 D$. Let $D_-$ be an open subset of $D$
which does not intersect a neighbourhood $U$ of $\6 D$,
and $D_+$ an open subset of $M \backslash D$ which
does not intersect $U$. Then there exists an LCK form
$\omega_1$ which is equal to $\omega$ on $D_+$
and to $dd^c \psi$ on $D_-$.

\hfill

\proof
We use the same proof as for \ref{_gluing_Kahler_Proposition_},
and note that the regularized maximum of 
$d_\theta d^c_\theta$-plurisubharmonic functions is
$d_\theta d^c_\theta$-plurisubharmonic by \ref{_d_theta_reg_max_Claim_}.
\endproof

\hfill

\remark\label{_non-strict_psh_Remark_}
\ref{_gluing_LCK_Proposition_} is  true also if  
$\psi$ and $\phi$ are not strictly $d_\theta
d^c_\theta$-pluri\-sub\-har\-monic.
In this case the gluing construction works,
but it gives a function which is $d_\theta
d^c_\theta$-plurisubharmonic, but not
 strictly $d_\theta d^c_\theta$-plurisubharmonic.


\section{Negative automorphic potentials for LCK metrics}\label{_nega_}


\theorem\label{_strictly_negative_pots_Theorem_}
Let $(M, \theta, \omega)$ be an LCK manifold
which is not K\"ahler,
and $\omega= d_\theta d^c_\theta(\phi)$ for some
smooth  $d_\theta d^c_\theta$-plurisubharmonic function $\phi$.
Then $\phi >0$ at some point of $M$.

\hfill

\proof Suppose, by absurd, that $\f\leq 0$
everywhere on $\tilde M$.
Since the LCK potential is stable under $C^2$-small deformations
of $\phi$, $\phi-\epsilon$ is also an LCK potential. Therefore, we may 
assume that $\phi <0$ everywhere. Define
$$\psi:=-\log(-\f).$$ 
Since $x\arrow -\log(-x)$ is strictly
monotonous and convex, the function
$\psi$ is strictly plurisubharmonic.  Moreover,
for every $\gamma\in\Gamma$, we have
\[
\gamma^*\psi=-\log(-(\f\circ\gamma))=-
\log(\chi(\gamma))-\log(-\f)=\const+\psi.
\]
Therefore, the K\"ahler form
$dd^c\psi$ is $\Gamma$-invariant and descends to $M$.
\endproof


\section{LCK potentials on Stein manifolds}
\label{_last_}


\subsection{Submanifolds with strictly pseudoconvex 
boundary and positivity of LCK potentials}

Let $(M, \theta, \omega)$ be an LCK manifold,
with $\omega= d_\theta d^c_\theta(\phi)$ for some
smooth $d_\theta d^c_\theta$-plurisubharmonic function $\phi$.
The condition $\omega= d_\theta d^c_\theta(\phi) >0$
is open in $C^2$-topology on the set of all functions $\phi$ on $M$. 
Adding a $C^2$-small function to $\phi$ if necessary,
we may assume that $0$ is a regular value of $\phi$.
The pullback of $\phi^{-1}(0)$ to $\tilde M$ is the set
of zeros of a K\"ahler potential. Therefore, 
it is strictly pseudoconvex, and the same is true
about  $\phi^{-1}(0)$. Since $\phi^{-1}(c)$ are $C^2$-close
to $\phi^{-1}(0)$ as subvarieties for small $c$, these sets are also
strictly pseudoconvex.

Choose a regular value $c >0$ of $\phi$ such that $\phi^{-1}(c)$
is non-empty and pseudoconvex. Then $\phi^{-1}(c)$ is a strictly pseudoconvex
CR-submanifold in $M$, and $D:= \phi^{-1}(]-\infty, c])$
is a strictly pseudoconvex set with boundary. Note that the interior of $D$ is an open submanifold in $M$, and hence it is LCK.

Then our main result (\ref{_main_intro_Theorem_})
follows from the gluing theorem \\ (\ref{_gluing_LCK_Proposition_})
and the following result about LCK manifolds with pseudoconvex
boundary.

\hfill

\theorem\label{_boundary_LCK_psh_Theorem_}
Let $(M,\theta, \omega)$ be a compact
LCK manifold of LCK rank 1 with smooth boundary which is strictly
pseudoconvex. Assume that $\omega=  d_\theta d^c_\theta(\phi)$
for some smooth $d_\theta d^c_\theta$-plurisubharmonic function $\phi$.
Then $M$ admits a positive 
$d_\theta d^c_\theta$-plurisubharmonic function $\psi$
such that $\psi$ is constant on the boundary $\6 M$.

\hfill

We prove \ref{_boundary_LCK_psh_Theorem_} later in this section.
Let us deduce \ref{_main_intro_Theorem_}
from \ref{_boundary_LCK_psh_Theorem_} and gluing.

\hfill

\theorem\label{_main_Theorem_}
Let $(M, \omega, \theta)$ be an LCK manifold manifold with a K\"ahler covering 
admitting an automorphic K\"ahler potential $\tilde \phi$.
Then $M$ also admits an LCK metric with a
positive automorphic potential and the same Bott-Chern class of the fundamental form.

\hfill

\proof Let $\phi$ be the corresponding potential on $M$,
$d_\theta d^c_\theta\phi=\omega$. Then $\phi>0$ somewhere
on $M$ (\ref{_strictly_negative_pots_Theorem_}).
As above, choose a regular value $c >0$ of $\phi$ such that $\phi^{-1}(c)$
is non-empty, and let $D:= \phi^{-1}(]-\infty, c])$
be the corresponding  strictly pseudoconvex set with boundary. 
By \ref{_boundary_LCK_psh_Theorem_},  $D$ admits a positive $d_\theta d^c_\theta$-plurisubharmonic function $\psi$
such that $\psi$ is constant on the boundary $\6 M$.
Choosing $\epsilon>0$ sufficiently small, and modifying
$\psi$ by adding a $C^2$-small function for transversality,
we may assume that the set
$S:=\{ m\in M \ \ |\ \ \epsilon \psi(m)=\phi(m) > 0\}$
is smooth and compact in $\phi^{-1}([0,c])$,
and $\Lie_X\phi > \epsilon \Lie_X \psi$ on $S$
as in \ref{_gluing_LCK_Proposition_}.
Then we may glue $\phi$ and $\epsilon \psi$
(\ref{_gluing_LCK_Proposition_}, \ref{_non-strict_psh_Remark_}).
We obtain an everywhere positive  $d_\theta
d^c_\theta$-plurisubharmonic 
function $\phi_1$. Adding $\delta \phi$, for $\delta >0$
sufficiently small, we make sure that  $\phi_1 + \delta
\phi$ is everywhere positive and strictly $d_\theta
d^c_\theta$-plurisubharmonic.
\endproof

\subsection{LCK potentials on 
submanifolds with strictly pseudoconvex boundary}

To finish the proof of the main theorem, it remains
to construct positive LCK potentials on LCK manifolds
with pseudoconvex boundary (\ref{_boundary_LCK_psh_Theorem_}).

\hfill

\theorem\label{_boundary_LCK_psh-bis_Theorem_}
Let $(D,\theta, \omega)$ be a compact 
LCK manifold of LCK rank 1 with smooth boundary which is strictly
pseudoconvex. Assume that $\omega=  d_\theta d^c_\theta(\phi)$
for some smooth $d_\theta d^c_\theta$-plurisubharmonic
function $\phi$, such that for some $\epsilon >0$
the function $\phi-\epsilon$ is also
$d_\theta d^c_\theta$-plurisubharmonic,
vanishes on the boundary of
$D$, and is strictly negative on $D$.
Then $D$ admits a positive
$d_\theta d^c_\theta$-plu\-ri\-sub\-har\-monic function $\psi$
which is constant on the boundary $\6 D$.

\hfill

{\bf Proof:}
Note that a manifold $D$ with smooth,
strictly pseudoconvex boundary is holomorphically convex. 
Then the Remmert reduction (\ref{rem_red}) implies that $D$
admits a proper, surjective  and holomorphic 
map $\pi:\; D \arrow D_0$ with connected fibres to a Stein
variety $D_0$ with isolated singularities.

Let $\Phi$ be a negative K\"ahler potential on the covering
$\tilde D$ of $D$, obtained from $\phi-\epsilon$.
Then the function $\Psi:=-\log(-\Phi)$ is strictly plurisubharmonic
on $\tilde D$ (see \ref{_strictly_negative_pots_Theorem_}),
hence the 1-form $\Theta:=d^c\Psi$ is defined on $D$.
On the other hand, $d\Theta$ is a K\"ahler form.
This implies that $D$ has no compact subvarieties (without boundary),
and the map $\pi:\; D \arrow D_0$ is bijective.
This implies that $D$ is Stein.

Now, let $L$ be the weight
bundle on $D$, associated to the character \eqref{chi}.
Denote by $\tilde D$ the smallest K\"ahler covering of
$D$. Then $\tilde D$ is a $\Z$-covering of $D$, and $L$ is
trivial on $\tilde D$. 
We call a function $f$ on $\tilde D$ {\bf
  $\lambda$-automorphic} if for each $\gamma\in \pi_1(D)$,
we have $\gamma^*(f)= \chi(\gamma)^\lambda f$.


Clearly,  
$1$-automorphic holomorphic functions on $\tilde D$ correspond to 
holomorphic sections of $L$.
Since $D$ is Stein,
the space of sections of a holomorphic bundle is 
globally generated; this assures the existence of
sufficiently many holomorphic sections of $L$.
Then, for a sufficiently big collection $f_i$ of sections
of $L$, the sum $\sum |f_i|$ is positive everywhere on
$D$. This gives an 1-automorphic plurisubharmonic
function on $\tilde D$.

We have proven that $D$ admits a positive
LCK potential. To obtain a potential which
is constant on $\6 D$, we perform the following trick.

Let $\{f_1, ..., f_n\}$ be the set of holomorphic
sections of $L$ without common zeros on $\6 D$.
Such a set exists and is nonempty because the pushforward of $L$
to a Stein variety $D_0$ is globally generated.
The following claim finishes the proof of
\ref{_boundary_LCK_psh-bis_Theorem_}, because
$|f_i| \cdot |b|$ is the length of 
the section $f_ib$ of $L$ for any holomorphic function $b$ on $D$.

\hfill

\claim \label{_sum_of_holo_Claim_}
Let $\{a_1, \ldots, a_n\}$ be a collection of non-negative
functions on a complex manifold $D$ with a smooth strictly
pseudoconvex boundary. Assume that $a_i$ have no common
zeroes on the boundary. Then for any positive function
$A$ on the boundary there exists a collection 
$\{b_1,\ldots, b_n\}$ of positive functions
such that $A=\sum a_i b_i$ and each $b_i$
can be obtained as the limit of a sum of absolute values
of holomorphic functions.

\hfill

{\bf Proof. Step 1:}
By a theorem of Bremmermann (\cite[Theorem 2]{_Bremmermann:equivalence_}),
every positive plurisubharmonic function on a pseudoconvex
manifold is Hartogs, that is, it belongs to the closure of the
cone generated by absolute values of holomorphic funcions.
Therefore, it would suffice to find a sum $A=\sum a_i b_i$
with $b_i$ positive, continuous and plurisubharmonic.

\hfil

{\bf Step 2:} By another theorem of Bremmermann
(\cite[Theorem 7.2]{_Bremmermann:Dirichlet_}), any function
on the boundary of a bounded holomorphically
convex domain can be extended inside to a plurisubharmonic
function $b$. Applying this result to $\log A$ and then
taking the exponential, we can make sure that $b$ is positive.
To prove \ref{_sum_of_holo_Claim_}, it remains
to find a collection of positive continuous functions
$\{b_1,\ldots, b_n\}$ on the boundary of $D$ such that 
$A=\sum a_i b_i$.

\hfill

{\bf Step 3:}
Since $a_i$ are non-negative and have no common zeros, 
their sum $B$ is positive. Then $\sum_i AB^{-1} a_i= AB^{-1}
\sum_i a_i =A$.
\endproof

\hfill

\noindent{\bf Acknowledgment:} We are grateful to Victor Vuletescu
for an interesting counterexample which stimulated our work on this
problem, Stefan Nemirovski for stimulating discussions and reference
to Bremermann, Cezar Joi\c ta for a careful reading of the paper,  Matei Toma  for communicating us the simple proof of \ref{_strictly_negative_pots_Theorem_}, and to Jason Starr
for invaluable answers given in Mathoverflow.

\hfill

{\small

\noindent {\sc Liviu Ornea\\
University of Bucharest, Faculty of Mathematics, \\14
Academiei str., 70109 Bucharest, Romania}, and:\\
{\sc Institute of Mathematics "Simion Stoilow" of the Romanian
Academy,\\
21, Calea Grivitei Str.
010702-Bucharest, Romania\\
\tt lornea@fmi.unibuc.ro, \ \  Liviu.Ornea@imar.ro}

\hfill

\noindent {\sc Misha Verbitsky\\
Laboratory of Algebraic Geometry, \\
Faculty of Mathematics, National Research University HSE,\\
7 Vavilova Str. Moscow, Russia}, also: \\
{\sc Universit\'e Libre de Bruxelles, D\'epartement de Math\'ematique\\
Campus de la Plaine, C.P. 218/01, Boulevard du Triomphe\\
B-1050 Brussels, Belgium\\
\tt verbit@verbit.ru, mverbits@ulb.ac.be }
}

{\scriptsize

\begin{thebibliography}{100}


\bibitem[Be]{bel} F.A. Belgun, {\em On the metric structure of non-K\"ahler complex surfaces}, Math. Ann. {\bf 317} (2000), 1--40.

\bibitem[Bre1]{_Bremmermann:equivalence_}
 H.J. Bremermann, {\em On the Conjecture of the Equivalence of the Plurisubharmonic Functions and the Hartogs Functions,}
Math. Ann. {\bf 131} (1956), 76--86.

\bibitem[Bre2]{_Bremmermann:Dirichlet_}
H.J. Bremermann, {\em On a generalized Dirichlet problem for plurisubharmonic
         functions and pseudoconvex domains. Characterization of Shilov boundaries,}   Trans. Amer. Math. Soc, {\bf 91} (1959), 246-276.

\bibitem[Bru]{bru} M. Brunella, {\em Locally conformally K\"ahler metrics on Kato surfaces}, Nagoya Math. J. {\bf 202} (2011), 77--81.  arXiv:1001.0530.

\bibitem[D1]{_Demailly_1982_}
J.-P. Demailly,
{\em Estimations $L^2$ pour l'op\'erateur $\bar\partial$ d'un fibr\'e vectoriel
holomorphe semi-positif au-dessus d'une vari\'et\'e
k\"ahl\'erienne compl\`ete,} Ann. Sci. Ecole Norm. Sup. 4e
S\'er. 15 (1982) 457--511.

\bibitem[DP]{_Demailly_Paun_}
Demailly, J.-P., Paun, M., {\em Numerical characterization of the
K\"ahler cone of a compact K\"ahler manifold}, Ann. of
Math. \textbf{159} (2004), 1247--1274. arXiv:math/0105176. 

\bibitem[DO]{_Dragomir_Ornea_}
  S. Dragomir, L. Ornea, Locally conformally K\"ahler manifolds, Progress in Math. {\bf 55}, Birkh\"auser, 1998.

\bibitem[GO]{go} P. Gauduchon, L. Ornea, {\em Locally conformally K\"ahler metrics on Hopf surfaces}, Ann. Inst. Fourier (Grenoble) {\bf 48} (1998), 1107--1128.

\bibitem[G]{goto} 
R. Goto, \emph{On the stability of locally conformal
  Kaehler structures},   J. Math. Soc. Japan {\bf 66} (2014), 
no. 4, 1375--1401.  
arXiv:1012.2285

\bibitem[GW]{_Green_Wu:approx_}
R. E. Greene and H. Wu, {\em $C^\infty$-approximations of convex,
subharmonic, and plurisubharmonic functions},
Ann. Scient. Ec. Norm. Sup. {\bf 12} (1979), 47--84.

\bibitem[MD]{ma_di} G. Marinescu, T. Dinh, {\em On the compactification of hyperconcave ends and the theorems of Siu-Yau and Nadel},  Invent. Math. {\bf 164}, No. 2,  (2006), 233--248. arXiv:math/0210485.  


\bibitem[OT]{ot} K. Oeljeklaus, M. Toma, {\em Non-K\"ahler compact complex manifolds associated to number  fields}, Ann. Inst. Fourier {\bf 55}, no. 1 (2005), 1291--1300.

\bibitem[OV1]{ov_ma_10} 
L. Ornea, M. Verbitsky, {\em Locally conformal K\"ahler manifolds with potential}, Math. Ann. {\bf 348} (2010), 25--33. arXiv:math/0407231

\bibitem[OV2]{ov_jgp_09}
L. Ornea, M. Verbitsky, {\em Morse-Novikov cohomology of locally conformally K\"ahler manifolds}, J. Geom. Phys. {\bf 59}, No. 3 (2009), 295--305. arXiv:0712.0107.

\bibitem[OV3]{_OV:_Topology_imrn_} 
L. Ornea and M. Verbitsky, 
\emph{Topology of Locally Conformally K\"ahler Manifolds
with Potential}, Int. Math. Res. Notices {\bf 4} (2010), 117--126.  arXiv:0904.3362.

\bibitem[OV4]{ov_imrn_12} L.  Ornea, M. Verbitsky, {\em Automorphisms of locally conformally Kaehler manifolds}, IMRN, vol. 2012 Nr. 4 894--903.  arXiv:0906.2836.

\bibitem[OV5]{_OV:emb_Sasakian_}
L.  Ornea, M. Verbitsky, {\em Embeddings of compact
  Sasakian manifolds},
Math. Res. Lett. 14 (2007), no. 4, 703--710.  arXiv:math/0609617.

\bibitem[OV6]{ov_jgp_16} 
L. Ornea, M. Verbitsky, {\em LCK rank of locally conformally K\"ahler manifolds with potential}, J. Geom. Phys. {\bf 107} (2016), 92--98. arXiv:1601.07413.

\bibitem[OV7]{ov_pams_16} L. Ornea, M. Verbitsky, {\em Locally conformally K\"ahler metrics obtained from pseudoconvex shells}, Proc. Amer. Math. Soc. {\bf 144} (2016), 325--335. arXiv:1210.2080.

\bibitem[OV8]{ov_hopf_surf} L. Ornea, M. Verbitsky, {\em
  Hopf surfaces in locally conformally Kahler manifolds
  with potential}, arxiv:1601.07421, v2, 10 pages.

\bibitem[OV9]{_OV:_embedding_new_proof_}
L. Ornea, M. Verbitsky,
{\em Embedding of LCK manifolds with potential into Hopf
  manifolds using Riesz-Schauder theorem}, 
 arXiv:1702.00985, 14 pages.

\bibitem[O]{O} 
A. Otiman, {\em Morse-Novikov cohomology of locally conformally K\"ahler surfaces},  arXiv:1609.07675.
  
\bibitem[Re]{rem} R. Remmert, {\em Sur les espaces analytiques holomorphiquement s\'eparables et holomorphiquement
convexes}, C. R. Acad. Sci. Paris {\bf 243} (1956), 118--121.


\bibitem[Ve]{_Verbitsky_vanishing_}
M. Verbitsky, {\em Theorems on the
vanishing of cohomology for
locally conformally hyper-K\"ahler manifolds}, Proc. Steklov Inst. Math.
{\bf 246} (2004) 54--78, arXiv:math/0302219.


\end{thebibliography}

}
\end{document}